\documentclass{article}
\usepackage{authblk}
\usepackage{amssymb}
\usepackage{amsthm}
\usepackage{amsmath}
\newtheorem{rem}{Remark}[section]

\newtheorem{proposition}{Proposition}[section]

\newtheorem{lemma}[proposition]{Lemma}%[section]

\newtheorem{theorem}[proposition]{Theorem}%[section]

\newtheorem{corollary}[proposition]{Corollary}%[section]

\def\p{\noindent{\bf Proof. }}
\def\q{\hspace*{\fill}$\Box$\medskip}
\DeclareMathOperator{\diam}{diam}

\begin{document}                % beginning of text.

\title{Analytic Connectivity of $k$-uniform hypergraphs \footnotemark[1]}

\author[1]{An Chang}
\author[2]{Joshua Cooper}
\author[3]{Wei Li}
\affil[1]{School of Computer and Information Science, Fujian Agriculture and Forestry University,  Fuzhou, Fujian, 350002, P.~R.~China}
\affil[2]{Department of Mathematics, University of South Carolina, Columbia, SC, 29208, USA}
\affil[3]{Center for Discrete Mathematics and Theoretical Computer Science, Fuzhou University, Fuzhou, Fujian, 350003, P.~R.~China}

\maketitle

\begin{abstract}
In this paper, we study the analytic connectivity of a $k$-uniform hypergraph $H$, denoted by $\alpha(H)$. In addition to computing the analytic connectivity of a complete $k$-graph, we present several bounds on analytic connectivity that relate it with other graph invariants, such as degree, vertex connectivity, diameter, and isoperimetric number.
\end{abstract}

%\keywords{H-eigenvalues;  analytic connectivity; k-uniform hypergraph; tensor; bound.}

\maketitle
\footnotetext[1]{The research is supported by the National Natural Science Foundation of China(No. 11331003); The Science and Technology Project of The Education Department in Fujian Province (JA13117).}

\section{Introduction}

Like graphs, hypergraphs have many applications in various fields \cite{berge}. However, a preponderance of interesting computational problems concerning hypergraphs are NP-hard. Spectral graph theory has played a significant role in remedying the apparent intractability of hard graph problems by providing approximation algorithms and iterative numerical methods. Hence, a search for an analogous spectral hypergraph theory has been become the focus of many researchers in recent years. As graphs are related to matrices, hypergraphs are related to tensors.

Let $m$ and $n$ be two positive integers. A {\bf tensor} $\mathcal{T}$ over the complex field $\mathbb{C}$ of order $m$ and dimension $n$ is a multidimensional array with entries $a_{i_1,i_2,\cdots,i_m}\in\mathbb{C}$, for each $(i_1,\cdots,i_m\in [n]=\{1,\cdots,n\})$. Particularly, if $m=1$, $\mathcal{T}$ is a vector in $\mathbb{C}^n$. If $m=2$, then $\mathcal{T}$ is a square matrix with $n^2$ elements. When $m\geq3$, $\mathcal{T}$ is a ``higher-order'' tensor. Furthermore, if its entries are invariant under any permutation of their indices, then a tensor $\mathcal{T}$ is said to be {\bf symmetric}.

Since the eigenvalues of higher-order tensors were independently proposed by Qi\cite{qi2005} and Lim \cite{lim2005}, numerous contributions to a framework for understanding the spectra of $k$-uniform hypergraphs via tensors have appeared. In 2009, Bul$\grave{o}$ and Pelillo \cite{pelillo2009} gave new bounds for the largest eigenvalue for the adjacency tensor of a uniform hypergraph with respect to its clique number. In 2011, Hu and Qi \cite{hu2012} proposed a definition for the Laplacian tensor of an even uniform hypergraph and analyzed its properties. In the following year, Cooper and Dutle \cite{cooper} presented some spectral results concerning hypergraphs that closely parallel those in the spectral theory of 2-graphs. In \cite{xie2013}, definitions for the Laplacian tensor and signless Laplacian tensor of a $k$-uniform hypergraph were proposed that extend the definition of Laplacian matrices and signless Laplacian matrices, and their H-eigenvalues and Z-eigenvalues were studied. In \cite{xie2012, pearson}, the authors investigated the H-eigenvalues and the Z-eigenvalues for the adjacency tensor of a $k$-uniform hypergraph. Furthermore, there is a rich and more general theory of eigenvalues for nonnegative
tensors \cite{kchang2008, kchang2011,friedland2013,liu,yang2010,zhang2012}. In \cite{cooper, kchang2008, friedland2013, yang2010}, a Perron-Frobenius theory for nonnegative tensors was established, supplying a fundamental tool for the study of hypergraph spectra.

The second smallest eigenvalue of the Laplacian matrix of a 2-graph $G$, denoted by $\lambda_2(G)$, is called the ``algebraic connectivity'' of $G$. This eigenvalue plays an important role in spectral graph theory. Using the definition of Laplacian tensor in \cite{xie2013}, Qi \cite{qi2013} defined a natural analytic connectivity of $k$-graphs. In this paper, we investigate upper and lower bounds on this parameter expressed in terms of the degree sequence, vertex connectivity, isoperimetric number and diameter. We also compute the the analytic connectivity of complete $k$-graphs.

The remainder of the present paper is organized as follows. In Section 2, some preliminary definitions concerning tensors and hypergraphs are given.  Moreover, we present several results concerning the algebraic connectivity of 2-graphs, which we compare to the analytic connectivity of $k$-graphs in section 3. We demonstrate in section 3 that these invariants have many similar properties, but differ in several respects.

 \section{Preliminaries}
Let $\mathcal{T}$ be a real tensor of order $m$ dimension $n$ and $x=(x_1,x_2,\cdots,x_n)$ be a vector in $\mathbb{C}^n$. Then
\begin{equation}\label{def}
  \mathcal{T}x^m:=\sum_{i_1,\cdots,i_m=1}^{n}a_{i_1,\cdots,i_m}x_{i_1}\cdots x_{i_m},
\end{equation}
 and $\mathcal{T}x^{m-1}$ is a vector in $\mathbb{C}^n$, whose $i$-th component is defined by
$$ (\mathcal{T}x^{m-1})_{i}=\sum_{i_2,\cdots,i_m=1}^{n}a_{i,i_2,\cdots,i_m}x_{i_2}\cdots x_{i_m},\;\;\;for \;\;\;i\in[n].$$
Let $x^{[r]}=(x_1^r, x_2^r, \ldots, x_n^r)$ be
a vector in $\mathbb{C}^n$, where $r$ is a positive integer.
If there is a non-zero vector $x\in \mathbb{C}^n$ and $\lambda\in \mathbb{C}$, satisfying
\begin{equation}\label{heigen}
\mathcal{T}x^{m-1}=\lambda x^{[m-1]},
\end{equation}
then one can say $\lambda$ is an {\bf eigenvalue} of $\mathcal{T}$, and $x$ is its corresponding {\bf eigenvector}.
In particular, if $x$ is real, then $\lambda$ is also real.
In this case, $\lambda$ is an $\mathbf{H-eigenvalue}$ of $\mathcal{T}$. And
 $x$ is its corresponding $\mathbf{H-eigenvector}$\cite{qi2004}. Let $\mathbb{R}^n_+$ and $\mathbb{R}^n_{++}$ be the set of all nonnegative vectors and the set of all positive vectors in $\mathbb{R}^n$, respectively. If $x\in \mathbb{R}^n_+$, then $\lambda$ is called an $\mathbf{H^+-eigenvalue}$ of $\mathcal{T}$. If $x\in \mathbb{R}^n_{++}$, then $\lambda$ is an
 $\mathbf{H^{++}-eigenvalue}$ of $\mathcal{T}$. If all the entries are nonnegative, then $\mathcal{T}$ is called a {\bf nonnegative tensor}.

In what follows, we employ standard definitions and notation from hypergraph theory; see, e.g.,\cite{berge}. A {\bf hypergraph} $H$ is a pair $(V,E)$. The elements of $V=V(H)=\{v_1,v_2,\cdots,v_n\}$ are referred to as {\bf vertices} and the elements of $E=E(H)=\{e_1,e_2,\cdots,e_m\}$ are called {\bf edges}, where $e_i\subseteq V$ for $i\in[m]$. A hypergraph $H$ is said to be {\bf $k$-uniform} for an integer $k\geq2$, if for all $e_i\in E(H)$, $|e_i|=k$, where $i\in[m]$. We often use the term {\bf $k$-graph} in  place of ``$k$-uniform hypergraph'' for short. Clearly, a $2$-graph is what is usually termed a simple, undirected graph. For any vertex $v_i\in V$ ($i\in[n]$), its degree $d(v_i)$ is defined as
$d(v_i)=|\{e_p:v_i\in e_p\in E\}|$. Denote the maximum degree, the minimum degree and the average degree
of $H$ by $\Delta(H)$, $\delta(H)$ and $\overline{d}(H)$, respectively. If $\overline{d}(H)=\Delta(H)$, then $H$ is a {\bf regular}
$k$-graph. Two vertices $v_i$ and $v_j$ are called {\bf adjacent} if and only if there exists an edge $e\in E(H)$ such that $\{v_i,v_j\}\subseteq e\in E(H)$. Two vertices $v_i$ and $v_j$ are called {\bf connected} if either $v_i$ and $v_j$ are adjacent, or there are vertices $v_{i_1},\cdots, v_{i_s}$ such that $v_i$ and $v_{i_1}$, $v_{i_s}$ and $v_j$, $v_{i_r}$ and $v_{i_{r+1}}$ for $r = 1,\cdots,s-1$ are adjacent, respectively. A $k$-graph $H$ is called {\bf connected} if any pair of its vertices are  connected. 

 The {\bf isoperimetric number} for a $k$-graph $H$, denoted by $i(H)$, is defined by
\begin{equation}\label{iso}
i(H)=\min \left \{\frac{|E_H(S,\overline{S})|}{|S|}: S\subset V(G), 0\leq S\leq \frac{|V(G)|}{2} \right \},
\end{equation}
where $\overline{S}=V\backslash S$ and the edge set $E_H(S,\overline{S})$ consist of the edges in $H$ with end vertices in both $S$ and $\overline{S}$. Particularly, when $k=2$, the edge $e\in E_H(S,\overline{S})$ has exactly one vertex in $S$ and one vertex in $\overline{S}$. When $k\geq3$, the edge $e\in E_H(S,\overline{S})$ satisfies $e\cap S\neq \emptyset$ and $e\cap \overline{S} \neq \emptyset$. For the sake of convenience,  $E_H(S)$,  or $E(S)$ for short, denotes the edge set  consisting of
edges in $H$ whose vertices are all in $S$.  Sometimes $E_H(S,\overline S)$ is called an {\bf edge cut} of $H$.  Indeed, if we delete $E(S,\overline S)$ from $H$, then $H$ is separated into two k-graphs $H[S] = (S,E(S))$ and $H[\overline S] = (\overline S,E(\overline S))$. The minimum cardinality of such an edge cut is called the {\bf edge connectivity} of $H$, denoted by $e(H)$. A {\bf vertex cut} of $H$ is a vertex subset $V'\subset V(H)$  such that $H-V'$ is disconnected, where $H-V'$ is the graph obtained by deleting all vertices in $V'$ and all incident edges. The {\bf vertex connectivity} of $H$, denoted by $v(H)$ is the minimum cardinality of any vertex cut $V'$. The complete $k$-graph has no vertex cut.

For a $k$-graph $H$ with $n$ vertices, the {\bf (normalized) adjacency tensor}, denoted by $\mathcal{A}(H)$ or $\mathcal{A}$ for short, is a tensor of order $k$ dimension $n$ with entries
\begin{equation}\label{adjamatrix}
  a_{i_1,i_2\cdots,i_k}=
 \left\{
  \begin{array}{ll}
  \frac{1}{(k-1)!} & if \;\;\{v_{i_1},v_{i_2},\cdots,v_{i_k}\}\in E(H)\\
  0 & otherwise.
  \end{array}
\right.
\end{equation}
The {\bf degree tensor} of $H$, denoted by $\mathcal{D}(H)$ or $\mathcal{D}$ for short, is a
diagonal tensor of order $k$ and dimension $n$ with its $i$-th diagonal entry as
$d_H(v_i)$ and $0$ otherwise. The {\bf Laplacian tensor}
and {\bf signless Laplacian tensor} of $H$ are defined as
$\mathcal{L}=\mathcal{L}(H)=\mathcal{D}(H)-\mathcal{A}(H)$ and
$\mathcal{Q}=\mathcal{Q}(H)=\mathcal{D}(H)+\mathcal{A}(H)$, respectively\cite{xie2013}.
It is easy to see that $\mathcal{A}$, $\mathcal{L}$ and
$\mathcal{Q}$ are all symmetric. Particularly, when $k=2$, the tensors
$\mathcal{A}$, $\mathcal{L}$ and $\mathcal{Q}$ are the adjacency
matrix $A$, the Laplacian matrix $L$ and the signless Laplacian
matrix $Q$ of a 2-graph, respectively.
% Therefore, the definition
%is simple and natural.

Let $e=\{v_{i_1},v_{i_2},\cdots,v_{i_k}\}\in E(H)$. Denote
\begin{equation*}
\mathcal{L}(e)x^k=\sum_{j=1}^{k}x_{i_j}^k-kx_{i_1}x_{i_2}\cdots x_{i_k}.
\end{equation*}
Then
\begin{equation*}\label{equ1}
\mathcal{L}x^k=\sum_{e\in E}\mathcal{L}(e)x^k.
\end{equation*}
Here we present the arithmetic-geometric mean inequality, which we refer to as the {\bf A-G inequality} for brevity.
\begin{lemma}
Let $a=(a_1,a_2,\cdots,a_n)$ be a vector in $\mathbb{R}^n_{+}$ and $A(a)=\frac{a_1+a_2+\cdots+a_n}{n}$, $G(a)=(a_1a_2\cdots a_n)^{\frac{1}{n}}$.
then $$A(a)\geq G(a).$$ Moreover, the equality holds if and only if $a_1=a_2=\cdots=a_n$.
\end{lemma}
 By A-G inequality, if $x\in \mathbb{R}^n_+$. then $\mathcal{L}(e)x^k\geq 0$ for each edge $e\in E(H)$, so that $\mathcal{L}x^k\geq 0$ as well.  Moreover, it is easy to verify that the smallest Laplacian H-eigenvalue of $H$ is exactly $0$ and the all-ones vector, denoted by $\mathbf{1}$, is the corresponding eigenvector. In \cite{qi2013}, Qi proved that
$$0=\min\{\mathcal{L}x^k: x\in \mathbb{R}_{+}^n,\sum_{i=1}^{n}x_i^k=1\}.$$
Qi also defined the {\bf analytic connectivity} $\alpha(H)$ of the $k$-graph $H$ by
\begin{equation}
\alpha(H)=\min_{j=1,2,\cdots,n}\min\{\mathcal{L}x^k:x\in \mathbb{R}_{+}^n,\sum_{i=1}^{n}x_i^k=1,x_j=0\}.
\end{equation}
Clearly, $\alpha(H)\geq 0$. The following statements illuminate the name of this parameter.
\begin{theorem}\cite{qi2013}\label{th3-1}
The $k$-graph $H$ is connected if and only if $\alpha(H)>0$.
\end{theorem}
\begin{theorem}\cite{qi2013}\label{3-1-1}
For a $k$-graph $H$, we have
$$e(G)\geq \frac{n}{k}\alpha(H).$$
\end{theorem}

Recall that the second smallest eigenvalue of Laplacian metric of a 2-graph $G$ of order $n$, denoted by $\lambda_2(G)$, is
defined by
\begin{equation}
\lambda_2(G)=\min_{\substack{x\perp \mathbf{1}\\ x\neq\mathbf{0}}}\frac{\left \langle Lx,x \right \rangle}{x^Tx},
\end{equation}
where $\mathbf{0}$ is the zero vector.
Fiedler obtained another important expression for $\lambda_2(G)$:
\begin{equation}\label{alge}
\lambda_2(G)=2n\min\left\{\frac{\sum_{v_iv_j\in E(G)}(x_i-x_j)^2}{\sum_{i=1}^{n}\sum_{j=1}^{n}(x_i-x_j)^2}: x\neq c\mathbf{1}, c\in \mathbb{R} \right\}.
\end{equation}
Fiedler referred to the number $\lambda_2(G)$ as the {\bf algebraic connectivity} of $G$, and related it to the classical connectivity of 2-graphs.
\begin{theorem}\label{th2-1}\cite{fiedler}
A 2-graph $G$ is connected if and only if $\lambda_2(G)>0$.
\end{theorem}
\begin{theorem}\label{2-1-1}\cite{fiedler}
For a 2-graph $G$, $\lambda_2(G)\leq v(G)\leq e(G)$.
\end{theorem}
There are several prominent graph classes for which the algebraic connectivity is known.  Here we give the algebraic connectivity when $G$ is a complete graph.
\begin{theorem}\label{th2-2}
If $K_n$ is a complete 2-graph of order $n$ with $\binom{n}{2}$ edges, then
\begin{align*}
\lambda_2(K_n)=n
\end{align*}
with corresponding eigenvector $x=(n,-1,-1,\cdots,-1)$.
\end{theorem}
There are also several bounds to the algebraic connectivity related to the parameters of a graph, such as degree, isoperimetric number, and diameter.
\begin{theorem}\label{th2-5}
Let $G$ be a $2$-graph with more than one edge and $(d(v_1),d(v_2),$ $\cdots,d(v_n))$ be the degree sequence of $G$. Then
\begin{equation}\label{dg2}
\lambda_2(G)\leq \min_{\{v_{i},v_{j}\}\in E(G)}\left\{\frac{d(v_{i})+d(v_{j})-2}{2}\right\}.
\end{equation}
\end{theorem}

 In general, the isoperimetric number is very hard to compute, and even obtaining any lower bounds on $i(H)$ seems to be a difficult problem. However, for 2-graphs $G$, the algebraic connectivity provides a reasonably good bound on $i(G)$. The following is a well-known inequality often called the ``{\bf Cheeger inequaility}''.
\begin{theorem}\cite{mohar}\label{th2-3}
Let $G$ be a $2$-graph. Then
$$2i(G) \geq \alpha(G)\geq \Delta(G)-\sqrt{\Delta(G)^2-i^2(G)}.$$
\end{theorem}
Regarding the {\bf diameter} of a k-graph $H$, denoted by $\diam(H)$, the following is a lower bound for 2-graphs.
\begin{theorem}\cite{mohar2}\label{th2-4}
Let $G$ be a 2-graph with $n$ vertices. Then $$\lambda_2(G)\geq \frac{4}{\diam(G)\cdot n }.$$
\end{theorem}

\section{Main results}
Evidently, the analytic connectivity and the algebraic connectivity are closely related to the connectivity of a graph. In this section, we will study the properties of analytic connectivity for $k$-graphs compared with that of algebraic connectivity for 2-graphs.

Given a finite set $X$ and positive integers $k$, $r$, and $\lambda$, a {\bf 2-design} (or {\bf balanced incomplete block design}) is a multiset of $k$-element subsets of $X$, called {\bf blocks}, such that the number of blocks containing any element of $X$ is $r$ and the number of blocks containing any pair of distinct $x,y \in X$ is $\lambda$.  If none of the elements of the multiset are repeated, then the 2-design is said to be {\bf simple}.  It is easy to see that a simple 2-design on a set of size $n$ with parameters $k$, $r$, and $\lambda$ as above is the same as a $k$-uniform hypergraph on $n$ vertices all of whose vertices have degree $k$, and all of whose codegrees $c(x,y) = |\{e \in E : x,y \in e\}|$ for $x \neq y$ are $\lambda$.

\begin{theorem}\label{th3-2}
	Let $H$ be a simple $2$-design with $n$ vertices, $b$ blocks, $k$ vertices per block, $r$ blocks containing each vertex, and $\lambda$ blocks containing each pair of vertices.  Further suppose that $H$ has no cut-vertex.  Then $$\alpha(H)=\lambda,$$
	with corresponding vector $x=\left((\frac{1}{n-1})^{\frac{1}{k}},(\frac{1}{n-1})^{\frac{1}{k}},\cdots,(\frac{1}{n-1})^{\frac{1}{k}},0\right)$.
\end{theorem}
\p It is well known that $bk=nr$ and $\lambda(n-1)=r(k-1)$.  Let $x=(x_1,x_2,\cdots,x_n)$ be a vector satisfying
\begin{equation*}
x_i=\left\{
\begin{array}{ll}
(\frac{1}{n-1})^{\frac{1}{k}},& \textrm{ if }\;\;\;i=1,2,\cdots,n-1;\\
0,& \textrm{ if }\;\;\;i=n.
\end{array}
\right.
\end{equation*}
Then, for any edge $e\subseteq \{v_1,v_2,\cdots,v_{n-1}\}$, we have $\mathcal{L}(e)x^k=0$. For those edges $e$ containing vertex $v_n$, we have $\mathcal{L}(e)x^k=\frac{k-1}{n-1}$. Since there are $r$ edges containing $v_n$, it follows that
\begin{equation}\label{eq4}
\alpha(H)\leq \mathcal{L}x^k=r \frac{k-1}{n-1}=\lambda.
\end{equation}
On the other hand, suppose that $y=(y_1,y_2,\cdots,y_{n-1},y_n)\in \mathbb{R}^n_+$ is the vector achieving $\alpha(H)$; we may assume without loss of generality that $y_n=0$. According to the A-G inequality, for each edge $e\in E(H)$, $\mathcal{L}(e)y^k\geq 0$. In addition, for those edges $e$ containing the vertex $v_n$, we have $\mathcal{L}(e)y^k=\sum_{j=1}^{k-1}y_{i_j}^k$, where $i_j\neq n$.
Therefore,
\begin{equation}\label{eq3}
\alpha(H)\geq \lambda \sum_{i=1}^{n-1}y_i^k=\lambda.
\end{equation}
Combining (\ref{eq4}) and (\ref{eq3}), we have
\begin{equation*}
\alpha(H)=\lambda.
\end{equation*}
It remains to verify that any vector $x$ achieving $\alpha(H)$ has the desired form. If equality holds in (\ref{eq3}), then every edge $e\subseteq E(H)$ with $v_n \not \in e$ satisfies $\mathcal{L}(e)x^k=0$, since no edge containing $v_n$ contributes to the sum. By the A-G inequality, each coordinate of $x$ corresponding to a vertex in $e$ has the same value. Since the subgraph induced by $V(H)\setminus \{v_n\}$ is connected, we may conclude that $x_1=x_2=\cdots=x_{n-1}$. Moreover, $\sum_{i=1}^{n-1}x_i^k=1$, so that $x_1=x_2=\cdots=x_{n-1}=(\frac{1}{n-1})^{\frac{1}{k}}$.
The result follows. \q

We can now give the analytic connectivity of a complete $k$-graph $K_n^{(k)}$, as follows, since $K_n^{(k)}$ is a $2$-design with $\lambda = \binom{n-2}{k-2}$.

\begin{corollary}\label{cor3-2} $\alpha(K_n^{(k)})=\binom{n-2}{k-2}$.
\end{corollary}

\begin{rem}
We compare with Theorem \ref{th2-2} by taking $k=2$ in Theorem \ref{th3-2}, so that $\alpha(K_n)=1$.
\end{rem}

Two additional properties of $\alpha(H)$ are given in Theorem \ref{subgraphs1} and Corollary \ref{subgraphs}. Furthermore, Theorem \ref{cutset} presents an upper bound on $\alpha(H)$ in terms of vertex connectivity.
 
\begin{theorem}\label{subgraphs1}
If $H_1$ and $H_2$ are edge-disjoint hypergraphs with the same vertex set then $\alpha(H_1)+\alpha(H_2)\leq \alpha(H_1\cup H_2)$.
\end{theorem}

\p Since $H_1$ and $H_2$ are edge-disjoint, we have $\mathcal{L}(H_1\cup H_2)x^{k}=\mathcal{L}(H_1)x^k+\mathcal{L}(H_2)x^k$. Obviously, $\alpha(H_1\cup H_2)\geq \alpha(H_1)+\alpha(H_2)$.
\q

\begin{corollary}\label{subgraphs}
If $H_1$ and $H_2$ have the same vertex set and $E(H_1)\subseteq E(H_2)$, then $\alpha(H_1)\leq \alpha(H_2)$.
\end{corollary}

\begin{theorem}\label{cutset}
	Let $H$ be a hypergraph of order $n$ and $v(H)$ be the vertex connectivity of $H$. Then
	$$\alpha(H)\leq \binom{n-2}{k-2}-\left[\binom{n-v(H)-1}{k-1}-\binom{\frac{n-v(H)}{2}-1}{k-1}\right]\frac{k-1}{n-1}.$$
\end{theorem}

\p Let $S$ be a minimum cut set of $H$, {\it i.e.}, $|S|=v(H)$; and let $H_1,H_2,\cdots,H_l$ be the components of $H-S$, with $n_i = |V(H_i)|$ for $i\in [l]$. Without loss of generality, suppose $n_1 = \min_{i \in [l]} n_i$, so that $1\leq n_1\leq \frac{n-v(H)}{2}$ and there is a vertex $u\in V(H_1)$.

Define a vector $x=(x_1,x_2,\cdots,x_n)$ by
\begin{equation*}
x_i=\left\{
\begin{array}{ll}
\left(\frac{1}{n-1}\right)^{\frac{1}{k}} &\textrm{ if } \mbox{$v_i\neq u $;}\\
0 &\textrm{ if }\mbox{$v_i=u$.}
\end{array}
\right.
\end{equation*}
It is easy to see $\sum_{i=1}^{n}x_i^{k}=1$. Let $H'$ arise from $H$ by adding all edges $e$ containing $u$ such that $e \cap S \neq \emptyset$ or $e \cap H_i = \emptyset$ for all $i \in \{2,\ldots,l\}$; note that all edges of $H$ have this form, so $H \subset H'$.  The maximum possible number of edges containing $u$ is $n-1 \choose k-1$.  However, we exclude from $H'$ those edges $e$ such that $e \cap S= \emptyset$ and $e \cap (\bigcup_{j=2}^{l}V(H_j)) \neq \emptyset$. Therefore,  $d_{H'}(u)=\binom{n-1}{k-1}-\binom{n-v(H)-1}{k-1}+\binom{n_1-1}{k-1}$. Moreover, $\mathcal{L}(e)x^k=\frac{k-1}{n-1}$ if $e$ contains vertex $u$ and $\mathcal{L}(e)x^k=0$ otherwise. Then, by Corollary \ref{subgraphs}, we have
\begin{align*}
\alpha(H)\leq \alpha(H')\leq&\mathcal{L}(H')y^k\\
=&\left[\binom{n-1}{k-1}-\binom{n-v(H)-1}{k-1}+\binom{n_1-1}{k-1}\right]\frac{k-1}{n-1}\\
\leq& \binom{n-2}{k-2}-\left[\binom{n-v(H)-1}{k-1}-\binom{\frac{n-v(H)}{2}-1}{k-1}\right]\frac{k-1}{n-1}
\end{align*}
The result follows.
\q

\begin{rem}
Taking $k=2$, {\it i.e.,} for a 2-graph $G$, we have $\alpha(G) \leq \frac{n+v(G)-2}{2(n-1)}.$
\end{rem}

Next, we investigate upper and lower bounds on $\alpha(H)$ in terms of the isoperimetric number and diameter. Before coming to our results, two extended A-G inequalities are needed.
\begin{lemma}
Let $a=(a_1,a_2,\cdots,a_n)$ be a vector in $\mathbb{R}^n_{+}$ and $A(a)=\frac{a_1+a_2+\cdots+a_n}{n}$, $G(a)=(a_1a_2\cdots a_n)^{\frac{1}{n}}$.	
\begin{item}
\item{(1)} \cite{kuang} Suppose that $b_j=a_{\sigma(i)}$, where $j=1,2,\cdots,n$ and $\sigma\in S_n$ is a permutation of the set $[n]$, then
\begin{equation}\label{ineq1}
A(a)-G(a)\geq \frac{1}{n}\sum_{j=1}^{m}(\sqrt{b_j}-\sqrt{b_{n+1-j}})^2,
\end{equation}
where $m=\lfloor\frac{n}{2}\rfloor$. Moreover, equality holds if and only if $b_1b_n=b_jb_{n+1-j}$, $1\leq j\leq m$.
\item{(2)}
\begin{equation}\label{ineq2}
A(a)-G(a)\geq \frac{1}{(n-1)n}\sum_{1\leq i<j\leq n}(\sqrt{a_i}-\sqrt{a_{j}})^2.
\end{equation}
Moreover, equality holds if and only if $a_1=a_2=\cdots=a_n$.
 \end{item}
\end{lemma}

\p Here, it is sufficient to verify (\ref{ineq2}). Since
$$A(a)=\frac{1}{(n-1)n}\sum_{1\leq i<j\leq n}(b_i+b_j),$$
$$G(a)=\left(\prod_{1\leq i<j\leq n}\sqrt{b_ib_j}\right)^{\frac{2}{(n-1)n}},$$
then,
\begin{align*}
&A(a)-\frac{1}{(n-1)n}\sum_{1\leq i<j\leq n}(\sqrt{b_i}-\sqrt{b_j})^2\\
=&\frac{2}{(n-1)n}\sum_{1\leq i<j\leq n}\sqrt{b_ib_j}\\
\geq&\left(\prod_{1\leq i<j\leq n}\sqrt{b_ib_j}\right)^{\frac{2}{(n-1)n}}=G(a).  \;\;\;\;\; \mbox{(by A-G inequality)}
\end{align*}
Hence, the inequality (\ref{ineq2}) holds. Moreover, equality holds if and only if $b_ib_j=b_sb_t$, for any $i,j,s,t\in[n]$, which means $b_1=b_2=\cdots=b_n$.\q

\begin{theorem}\label{main}
Let $H$ be a $k$-graph, where $k\geq3$. Then
$$\frac{k}{2} i(H) \geq \alpha(H)\geq \Delta-\sqrt{\Delta^2-i^2(H)}.$$
\end{theorem}
\p Suppose that the isoperimetric number $i(H)$ is witnessed by the set $S$. Let $y=(y_1,y_2,\cdots,y_n)$ be a vector in $\mathbb{R}^n_+$ satisfying
\begin{equation*}
y_i=\left\{
\begin{array}{ll}
1/(|S|)^{\frac{1}{k}}& v_i\in S\\
0& v_i\in \overline{S}.
\end{array}
\right.
\end{equation*}
Denote $t(e)=|\{v|v\in e \cap S\}|$ and $t(S)=\frac{\sum_{e\in E(S,\overline{S})}t(e)}{|E(S,\overline{S})|}$.
Then
\begin{align}\label{eq1}
\alpha(H)\leq \mathcal{L}y^k=& \left (\sum_{e_p\in E(S)}+\sum_{e_p\in E(\overline{S})}+\sum_{e_p\in E(S,\overline{S})} \right)\mathcal{L}(e_p)y^k\nonumber \\
=&\sum_{e_p\in E(S,\overline{S})}\sum_{v_i\in e_p\cap S}x_i^k\nonumber\\
\leq&t(S)\frac{1}{|S|}|E(S,\overline{S})|\nonumber\\
=&t(S)i(H).
\end{align}
Similarly, we have
\begin{equation}\label{eq2}
\alpha(H)\leq t(\overline{S})i(H).
\end{equation}
Summing (\ref{eq1}) and (\ref{eq2}), since $t(S)+t(\overline S)=k$, we have
$$\alpha(H)\leq \frac{k}{2}i(H).$$
Thus, the proofs for the Cheeger inequality regarding an upper bound for $\alpha(H)$ is completed.

On the other hand, to verify the lower bound of $\alpha(H)$ for a bipartition $(S, \bar{S})$ on $V(H)$, we first define a multiple 2-graph with possible loops, $\widehat{H}=(V(\widehat{H}),E(\widehat{H}))$, where $V(\widehat{H})=V(H)$ and the edge set of $\widehat{H}$ is derived in the following way. Suppose $x=\{x_1,x_2,\cdots,x_n\}$ is the vector achieving $\alpha(H)$. For each edge  $e=\{v_{i_1},v_{i_2},\cdots,v_{i_k}\}\in E(H)$, with loss of generality, let $x_{i_1}\geq x_{i_2}\geq \cdots \geq x_{i_k}$. Then  $E(\widehat{H})=\bigcup_{e\in E(H)}\{v_{i_j}v_{k+1-i_{j}}: j=1,2,\cdots,\lfloor\frac{k}{2}\rfloor\}$.

Since $x_i\geq 0$ for $i\in [k]$, then
\begin{align*}
\alpha(H)=&\sum_{e=\{v_{i_1},\cdots,v_{i_k}\}\in E}\left(\sum_{j=1}^{n}x_{i_j}^k-k x_{i_1}\cdots x_{i_k}\right)\\
\geq& \sum_{\substack {e=\{v_{i_1},\cdots,v_{i_k}\}\in E \\ x_{i_1}\geq x_{i_2}\geq \cdots\geq x_{i_k}}}\sum_{j=1}^{k}\left (\sqrt{x_{i_j}^k}-\sqrt{x_{k+1-i_{j}}^k}\right )^2 \;\;\;\;\mbox{(By (\ref{ineq1}))}\\
=& \sum_{v_i,v_j\in E(\widehat{H})} \left (\sqrt{x_{i}^k}-\sqrt{x_{j}^k} \right)^2\\
\end{align*}

Let $M=\sum_{(v_i,v_j)\in E(\widehat{H})} (x_{i}^{k/2}-x_{j}^{k/2})^2$ and $y=x^{[\frac{k}{2}]}$. Moreover, for the sake of convenience, we denote by $E'$ the set $E(\widehat{H})$. According to the proof of Theorem \ref{th2-3}, we have

\begin{align}\label{eq2-5}
M=&\frac{\sum_{v_i,v_j\in E'}(y_i-y_j)^2\sum_{v_i,v_j\in E'}(y_i+y_j)^2}{\sum_{v_i,v_j\in E'}(y_i+y_j)^2}\\\nonumber
\geq& \frac{(\sum_{v_i,v_j\in E'}|y_i^2-y_j^2|)^2}{\sum_{v_i,v_j\in E'}(y_i+y_j)^2} \,\,\,\text{(by Cauchy-Schwarz)}
\end{align}

Let $0=t_0<t_1<t_2<\cdots<t_h$ be all distinct values of $y_i$, $i=1,2,\cdots,n$. For $s=0,1,2,\cdots,h$, let $V_s=\{v_i\in V: y_i\geq t_s\}$. For each edge $e\in E_H(V_s,\overline{V_s})$, let  $\delta_s(e)=\min\{|V_s\cap e|,| \overline {V_s} \cap e |\}$. Denote $\delta(V_s)=\min\{\delta_s(e):e\in E_H(V_s,\overline {V_s})\}$ and $\delta(H)=\min_{s\in [h]}\{\delta(V_s)\}$. Then

\begin{equation}\label{eq2-6}
\begin{split}
\sum_{v_i,v_j\in E'}|y_i^2-y_j^2|=&\sum_{i=1}^{h}\sum_{\substack{v_iv_j\in E'\\ y_i\geq y_j}}(y_i^2-y_j^2)\\
=&\sum_{i=1}^{h}\sum_{\substack{y_i=t_r\\y_j=t_l,l<r}}(t_r^2-t_{r-1}^2+t_{r-1}^2-\cdots-t_{l+1}^2+t_{l+1}^2-t_l^2)\\
=&\sum_{i=1}^{h}\sum_{v_i\in V_s}\sum_{v_j\notin V_s}(t_s^2-t_{s-1}^2)\\
\geq&\sum_{i=1}^{h}\delta(V_s)E_H(V_s,\overline{V_s})(t_s^2-t_{s-1}^2)\\
\geq&\delta(H)i(H)\sum_{i=1}^{h}|V_s|(t_s^2-t_{s-1}^2)\\
=&\delta(H)i(H)\sum_{i=1}^{h}t_k^2(|V_s|-|V_{s+1}|)\\
=&\delta(H)i(H)\sum_{i=1}^{h}y_i^2.
\end{split}
\end{equation}

On the other hand,
\begin{equation}\label{eq2-7}
\begin{split}
 \sum_{v_iv_j\in E'}(y_i+y_j)^2=&2\sum_{v_iv_j\in E'}(y_i^2+y_j^2)-\sum_{v_iv_j\in E'}(y_i-y_j)^2\\
\leq &2\sum_{i=1}^{n}d_iy_i^2-\sum_{v_iv_j\in E'}(y_i-y_j)^2\\
\leq & 2\Delta(\widehat{H})\sum_{i=1}^{n}y_i^2-\sum_{v_iv_j\in E'}(y_i-y_j)^2 \\
=&(2\Delta(\widehat{H})-M)\sum_{i=1}^{n}y_i^2\\
\leq& (2\Delta(H)-M)\sum_{i=1}^{n}y_i^2
\end{split}
\end{equation}
where $\Delta$ denotes maximum degree. Combining (\ref{eq2-5}), (\ref{eq2-6}) and (\ref{eq2-7}), we obtain
$$M\geq\frac{\delta(H)^2i(H)^2}{2\Delta(H)-M}.$$ Therefore,
$$M\geq \Delta(H)-\sqrt{\Delta(H)^2-\delta(H)^2i(H)^2}.$$
Since $\delta(H)\geq 1$, we have $\alpha(H)\geq\Delta(H)-\sqrt{\Delta(H)^2-i(H)^2}$.\q

\begin{theorem}
Let $H$ be a $k$-graph. Then
\begin{equation}
\alpha(H)\geq \frac{4}{n^2(k-1)\diam(H)}.
\end{equation}
\end{theorem}

\p Let $x=(x_1,x_2,\cdots,x_n)$ be the vector achieving $\alpha(H)$, where $x_n=0$, and $y=x^{[\frac{k}{2}]}$. Define a multiple 2-graph $H^*$ as follows. It has vertex set $V(H)$, and vertices $u$ and $v$ are adjacent in $H^*$ if and only if $\{u,v\}\subset e\in E(H)$. Evidently, $\diam(H)=\diam(H^*)$. Therefore,

\begin{align*}
\alpha(H)=&\mathcal{L}x^k=\sum_{e\in E(H)}\mathcal{L}(e)x^k\\
=&\sum_{e=\{v_{i_1},v_{i_2},\cdots,v_{i_k}\}\in E(H)}\left(\sum_{j=1}^{k}x_{i_j}^k-kx_{i_1}x_{i_2}\cdots x_{i_k}\right)\\
\geq&\sum_{e=\{v_{i_1},v_{i_2},\cdots,v_{i_k}\}\in E(H)}\frac{1}{k-1}\sum_{1\leq s<t\leq k}\left (x_{i_s}^{k/2}-x_{i_t}^{k/2} \right)^2\;\;\;\mbox{(by (\ref{ineq2}))}\\
=&\frac{1}{k-1}\sum_{v_iv_j\in E(H^*)}\left (x_{i}^{k/2}-x_{j}^{k/2} \right)^2\\
=&\frac{1}{k-1}\sum_{v_iv_j\in E(H^*)}(y_{i}-y_{j})^2\\
=&\frac{1}{k-1}\sum_{i=1}^{n}\sum_{j=1}^{n}(y_i-y_j)^2
\frac{\sum_{v_iv_j\in E(H^*)}(y_i-y_j)^2}{\sum_{i=1}^{n}\sum_{j=1}^{n}(y_i-y_j)^2}\\
\geq&\frac{\lambda_2(H^*)}{2n(k-1)}\sum_{i=1}^{n}\sum_{j=1}^{n}(y_i-y_j)^2\;\;\;\;\mbox{(by (\ref{alge}))}
\end{align*}
And
\begin{align}
\sum_{i=1}^{n}\sum_{j=1}^{n}(y_i-y_j)^2=&\sum_{i=1}^{n}\sum_{j=1}^{n}y_i^2+\sum_{i=1}^{n}\sum_{j=1}^{n}y_j^2-2\sum_{i=1}^{n}\sum_{j=1}^{n}y_iy_j\nonumber\\
=&2n\left \langle y,y\right \rangle^2-2\left \langle y,\mathbf{1}\right \rangle^2\nonumber\\
=&2n-2\left (\sum_{i=1}^{n-1}y_i \right)^2\nonumber\\
\geq&2n-2(n-1)\sum_{i=1}^{n-1}y_i^2 \nonumber\;\;\;\;\mbox{(by Cauchy-Schwarz)}\\
=&2.
\end{align}
Hence, from Theroem \ref{th2-4},
\begin{align*}
\alpha(H)\geq \frac{\lambda_2(H^*)}{n(k-1)}\geq \frac{4}{n^2(k-1)\diam(H^*)}=\frac{4}{n^2(k-1)\diam(H)}.
\end{align*}
Completing the proof.\q

The last theorem gives an upper bound on the analytic connectivity of $k$-graphs as a function of degree sequence.

\begin{theorem}
	Let $H$ be a $k$-graph with more than one edge. Then
	\begin{equation}\label{dg}
	\alpha(H)\leq \min\left\{\frac{d(v_{i_1})+d(v_{i_2})+\cdots+d(v_{i_k})-k}{k}:\{v_{i_1},v_{i_2},\cdots,v_{i_k}\}\in E(H)\right\}.
	\end{equation}
\end{theorem}
\p Let $e_p=\{v_{i_1},v_{i_2},\cdots,v_{i_k}\}$ be an edge in $H$ and $x=(x_1,x_2,\cdots,x_n)$ be a vector defined by
\begin{equation*}
x_i=\left\{
\begin{array}{ll}
k^{-1/k} & \textrm{ if } \;\;v_i\in e_p,\\
0& \textrm{ otherwise}.
\end{array}
\right.
\end{equation*}
Then $\sum_{i=1}^{n}x_{i}^k=1$. Moreover, it is easy to see that $\mathcal{L}(e_p)x^k=0$.
Then,
\begin{align*}
\alpha(H)\leq &\mathcal{L}x^k=\sum_{e\in E(H)}\mathcal{L}(e)x^k\\
=&(d(v_{i_1})-1)(1/k)+(d(v_{i_2})-1)(1/k)+\cdots+(d(v_{i_k})-1)(1/k)\\
=&\frac{d(v_{i_1})+d(v_{i_2})+\cdots+d(v_{i_k})-k}{k}
\end{align*}
completing the proof.\q
\begin{rem}
	When $k=2$, the upper bound in inequality (\ref{dg}) is
	$$\min\left\{\frac{d(v_{i_1})+d(v_{i_2})-2}{2}:\{v_{i_1},v_{i_2}\}\in E(H)\right\},$$
	which is exactly the upper bound in (\ref{dg2}).
\end{rem}


\begin{thebibliography}{99}

\bibitem{berge}
     C. Berge, Hypergraphs, North-Holland Mathematical Library 45, North-Holland, Amsterdam, 1989.
\bibitem{qi2004}
     L. Qi, Eigenvalues of a supersymmetric tensor and positive definiteness of an even degree multivariate form. Department of Applied Mathematics, The Hong Kong Polytechnic University, 2004.
\bibitem{pelillo2009}
S.R. Bul$\grave{o}$, M. Pelillo, New bounds on the clique number of graphs based on spectral hypergraph theory, in: T.St$\ddot{u}$tzle ed., Learning and Intelligent Optimization, Spring Verlag, Berlin, (2009) 45-48.
\bibitem{cooper}
     J. Cooper, A. Dutle, Spectra of uniform hypergraphs. Linear Algebra and its Application, 436 (2012) 3268-3292.

\bibitem{spectra}
     D. Cvetkovi$\acute{c}$, M. Doob, H. Sachs, Spectra of Graphs: Theory and Application, third ed., Johann Ambrosius Barth Verlag, 1995.
\bibitem{kchang2008}
    K.C. Chang, K. Pearson, T. Zhang, Perron Frobenius Theorem for nonnegative tensors, Communications in Mathematical Sciences, 6 (2008) 507-520.
\bibitem{kchang2011}
    K.C. Chang, K. Pearson, T. Zhang, Primitivity, the convergence of the NZQ method, and the largest eigenvalue for nonnegative tensors, SIAM J. Matrix Anal. Appl., 32 (2011) 806-819.
\bibitem{xie2013}
  J. Xie, A. Chang, On the H-eigenvalues of the singless Laplacian tensor for an even uniform hypergraph, Frontiers of Mathematics in China, 8 (2013) 107-128.
\bibitem{xieLAA}
  J. Xie, A. Chang, On the Z-eigenvalues of the adjacency tensors for uniform hypergraphs, Linear Algebra and its Application, 439 (2013) 2195-2204.
\bibitem{xie2012}
J. Xie, A. Chang, On the Z-eigenvalues of the signless Laplacian tensor for an even uniform hypergraph, Numerical Linear Algebra with Applications, 20 (2013) 1030-1045.
\bibitem{pearson}
K. Pearson and T. Zhang, On the spectral hypergraph theory of the adjacency tensor, Graphs and Combinatorics, (2013)DOI:10.1007/s00373-013-1340-x.
\bibitem{friedland2013}
S. Friedland, S. Gaubert, L. Han, Perron-Frobenius theorem for nonnegative multilinear forms and extensions, Linear Algebra and its Applications, 438 (2013) 738-749.
\bibitem{liu}
Y. Liu, G, Zhou, N.F. Ibrahim, An always convergent algorithm for the largest eigenvalue of an irreducible nonnegative tensor, Journal of Computational and Applied Mathematics, 235 (2010) 286-292.
\bibitem{hu2012}
S. Hu, L. Qi, Algebraic connectivity of an even uniform hypergraph, Journal of Combinational Optimization, 24 (2012) 564-579.
\bibitem{yang2010}
Q. Yang, Y. Yang, Further results for Perron-Frobenius Theorem for nonngeative tensors, SIAM J.Matrix Anal. Appl.,31 (2010) 2517-2530.
\bibitem{yang2011}
Q. Yang, Y. Yang, Further results for Perron-Frobenius Theorem for nonngeative tensors II, SIAM J.Matrix Anal. Appl.,32 (2011) 1236-1250.
\bibitem{shao}
   J. Shao, A general product of tensors with applications. Linear Algebra Appl.,439 (2013) 2350-2366.
\bibitem{lim2005}
L.H. Lim, Singular values and eigenvalues of tensors, a variational approach, in Proceedings 1st IEEE international workshop on computational advances of multitensor adaptive processing, (2005) 129-132.
\bibitem{qi2005}
   L. Qi, Eigenvalues of a real supersymmetric tensor, Journal of Symbolic Computiation 40 (2005) 1320-1324.
\bibitem{qi2013}
   L. Qi, Symmetric nonnegative tensors and copositive tensors, Linear Algebra and Its Applications, 439 (2013) 228-238.
\bibitem{qi2012}
   L. Qi, $H^+$-eigenvalues of Laplacian and signless Laplacian tensors, Department of Applied Mathematics, The Hong Kong Polytechnic University, July, 2013.
\bibitem{zhang2013}
L. Zhang, L. Qi, Z. Luo, The dominant eigenvalue of an essentially nonnegative tensor,Numerical Linear Algebra with Applications, 20 (2013) 929-941.
\bibitem{zhang2012}
L. Zhang, L. Qi, Y. Xu, Linear convergence of the LZI algorithm for weakly positive tensors, Journal of Computational Mathematics, 30 (2012) 24-33.
\bibitem{mohar}
B. Mohar, Isoperimetric numbers of graphs, J. Combin. Theory, Ser. B, 47 (1989).
\bibitem{mohar2}
B. Mohar, Eigenvalues, diameter, and mean distance in graphs, Graphs and Combinatories, 7 (1991) 53-64.
\bibitem{kuang}
J.C. Kuang, Applied inequality, 3rd Ed., Shangdong Sicence and Technology Press, China, (2004) 36 (in Chinese).
\bibitem{fiedler}
M. Fiedler, Algebraic connectivity of graphs, Czechoslovak Mathematical Journal, 23(98) 1973.
\end{thebibliography}
\end{document}